\documentclass[11pt]{article}
\usepackage[final]{epsfig}
\usepackage[left=3cm,right=3cm, top=2.5cm,bottom=2.5cm,bindingoffset=0cm]{geometry}
\usepackage{graphics}
\usepackage{amsmath}
\usepackage{amsfonts}
\usepackage{latexsym}
\usepackage{amssymb, mathabx}
\usepackage{amsthm}
\usepackage{graphicx}
\usepackage{epstopdf}
\usepackage{multicol,multirow}
\usepackage{hyperref, enumitem}
\usepackage{mathtools}

 \usepackage{relsize}

\usepackage[nottoc]{tocbibind}

\usepackage{amsthm}

\mathtoolsset{showonlyrefs}

\usepackage{tikz}
\usetikzlibrary{positioning}
\usetikzlibrary{decorations.text}
\usetikzlibrary{decorations.pathmorphing}
\usetikzlibrary{decorations.markings}

\usetikzlibrary{arrows} 
\tikzset{
    >=stealth',
    pil/.style={
           ->,
           thick,
           shorten <=2pt,
           shorten >=2pt,}
}
\tikzset{->-/.style={decoration={
  markings,
  mark=at position .7 with {\arrow{>}}},postaction={decorate}}}
  \tikzset{a/.style={decoration={
  markings,
  mark=at position .52 with {\arrow{angle 90}}},postaction={decorate}}}
\tikzset{-<-/.style={decoration={
  markings,
  mark=at position .4 with {\arrow{<}}},postaction={decorate}}}

\makeatletter
\def\thmhead@plain#1#2#3{%
  \thmname{#1}\thmnumber{\@ifnotempty{#1}{ }\@upn{#2}}%
  \thmnote{ {\the\thm@notefont#3}}}
\let\thmhead\thmhead@plain

\makeatother

\newcounter{AppCounter}

\def\restrict#1{\raise-.5ex\hbox{\ensuremath|}_{#1}}

\newtheorem{remark-definition}[lemma]{Remark-Definition}

\newtheorem{proposition-conjecture}[lemma]{Proposition-conjecture}

\theoremstyle{definition}

\newcommand{\proofend}{\hfill$\Box$\bigskip}

\newcommand{\low}[1]{\raise-.0ex\hbox{$\scriptstyle #1$}}
\newcommand{\high}[1]{\raise.5ex\hbox{$\scriptstyle #1$}}

\newcommand{\marginnote}[1]
{%\mbox{}\marginpar{\center{\hspace{0pt}\tiny{\bf#1}}}
}

\newcounter{ai}

\newcounter{bk}

\title {The Euler non-mixing made easy}

\author{Boris Khesin\thanks{ Department of Mathematics,
University of Toronto, Toronto, Canada,
\tt{khesin@math.toronto.edu}}
\\
}
\date{}

\begin{document}

\maketitle
\begin{abstract}
The non-transitivity without extra constraints in the Euler equation   in any dimension is almost evident and can be derived, e.g. from Morse theory.
\end{abstract}

%\tableofcontents

The classical Euler equation describes the motion of an invscid incompressible fluid filling a manifold $M$ as an evolution of its divergence-free velocity field $v$:
$$
\partial_t v+\nabla_v v=-\nabla p.
$$
Here $p$ is the pressure function determined by the equation itself along with the divergence-free condition  ${\rm div }\,v =0$. In this note we are particularly interested in the 
three-dimensional setting, ${\rm dim }\,M=3$, while the results are also extended to any dimension. One of the main problems of hydrodynamics is the description of properties of the dynamical system defined by the Euler equation in the appropriate space of  velocity fields. 
It is known that this is a Hamiltonian system with the Hamiltonian function given by the $L^2$-energy of the fluid, with the short-time existence for the corresponding flow for a sufficiently smooth initial $v$ (for $v$ in $C^k$ with non-integer $k>1$). 
While for finite-dimensional Hamiltonian systems on compact manifolds one always has  the Poincar\'e recurrence, 
for the 2D Euler on an annulus $M=S^1\times \mathbb R$ there are wandering solutions \cite{DEJ, Nad}, i.e. a  neighborhood in the space of initial conditions, such that solutions starting in that neighborhood will never return to it after some time. Note that the existence of wandering solutions in 3D  Euler, as well as in 2D Euler on an arbitrary $M$ without boundary, is still an open question.

 There is another related property of dynamics, transitivity. Namely, while wandering solutions emanating from a certain neighborhood never return to it after some time, non-transitivity of a dynamical system means that there are two different neighborhoods such that solutions 
from one of them never pass through the other. The mixing property is an even stronger property than transitivity, asking any neighbourhood not only to occasionally overlap with any other, but to always admit a nonzero intersection after some time. 
It turns out that the Euler equation does possess the properties of non-transitivity and non-mixing:
\medskip

{\bf Theorem 1.} {\it $(i)$ The 3D Euler equation on a compact $M$ 
is non-transitive and hence non-mixing: there are two open neighbourhoods 
in the $C^k, \, k\ge 1$ phase space of velocities, so that the Euler flow image of one of them will never intersect the other (as long as the flow exists). Such neighborhoods can be chosen within (null-homologous) divergence-free fields with any  helicity and  sufficiently high energy.

$(ii)$ The Euler equation on  a compact  $M$ of any dimension is non-transitive and hence non-mixing  in the $C^k, \, k\ge 1$ phase space of velocities (as long as the flow exists). 
}
\medskip

This property is based on the existence of various first integrals, and, in particular, on vorticity transport, one of remarkable properties of the  Euler equation: in 3D the vorticity field $w={\rm curl }\,v$ is frozen into the flow. 
A unifying idea for proving non-mixing in the 3D Euler equation in \cite{CT, KKP1, KKP2} was as follows:
find two neighbourhoods in the space of all velocity fields with some incompatible
topological properties of their vorticities, so that the Euler solutions
with initial conditions in one of them, while preserving this
property, would not be able to enter the other one. 
In \cite{KKP1} those neighbourhoods contained the fields whose vorticities
have many invariant tori in different isotopy classes (and then one was
applying KAM, which required high regularity, $k>4$). In \cite{CT} this was the property of vorticity
to be of contact type or not, which allowed lower smoothness ($k\ge 1$).  
The use of the number of hyperbolic zeros of vorticity to prove the Euler non-mixing, which is the key idea for establishing $(i)$,  was already suggested in \cite{KKP1}. 
As was remarked there, at that time it was unclear how to develop it in order to construct exact divergence-free vector fields 
with prescribed values of energy and helicity. Below we show that this can be easily achieved by utilizing ``vortex plugs", while preserving even more subtle continuous invariants, such as multiplicators of zeros of the vorticity field  in the $C^2$-setting.

\medskip

There are variations of the above formulation, see  \cite{CT, KKP1,KKP2}: a)~different smoothness, b)~specifying bounds for helicity and energy, c)~existence of a countably many neighbourhoods, d)~existence in a given homotopy class for nowhere vanishing vorticity fields, e)~local non-mixing close to steady solutions (the latter was the initial motivation for the non-mixing study:  prove that some solutions will never get close to steady ones).
\medskip

\proof
 $(i)$ We are proving  non-transitivity, from which non-mixing follows. We start with the $C^2$ case for velocity in 3D.
Let $w_0$ be an initial vorticity $C^1$ field. Assume that it has only non-degenerate zeros (and possibly nondegenerate periodic trajectories) in $M$ (and hence only a finite number of them).
Then there is a $ C^1$-small neighbourhood $U(w_0)$ of $w_0$, such that all fields from $U(w_0)$  have the same number of zeros and they all are non-degenerate.
Similarly, if another initial vorticity field $w_1$ has a different number of non-degenerate zeros, there is a small neighborhood $V(w_1)$ of fields with the same property and the flow $\phi^t(V(w_1))$ will never intersect $U(w_0)$. Thus the Euler flow is non-transitive. 
Note that the argument above requires  $C^1$-closeness for the vorticity field 
(and hence $C^2$ for the velocity field), which is the optimal smoothness for that Morse-type argument 
to distinguish between different number of zeros.
\smallskip

To lower the smoothness to $C^0$-regularity  for  vorticity (and hence $C^1$ for  velocity)  one  
compares vorticity without zeros with vorticity having nondegenerate zeros. 
Namely,  $C^0$-small perturbations of vorticities with nondegenerate zeros have {\it at least} as 
many zeros as the unperturbed ones, while vorticities $C^0$-close to
the ones without zeros will also have no zeros.  (Note that any three-dimensional $M$ admits 
a nonzero vorticity field.) Therefore the above argument still works for $C^0$-closeness for
vorticities,  thus giving the optimal smoothness for non-transitivity.
Namely, a non-degenerate zero of a vector field always has index $\pm 1$, and according to 
 the index theorem (which is essentially the local Intermediate Value Theorem 
 in the vector-function setting), it must persist for $C^0$-close perturbations (and globally sums to the Euler characteristic). 
Note that one can weaken the nondegeneracy assumption on the vorticity to just having a certain number 
of zeros of index $\pm1$. 
\smallskip

 Finally, by using the local insertion  of ``vortex plugs" into a given divergence-free field (to described in Example 2 below) 
one can generate new vorticity fields with pairs of new nondegenerate zeros and with an arbitrary helicity. 
This is a local construction which can be thought of as an insertion of a small rotor providing one-directional fast linking of trajectories inside a small (and isolated from everything else) invariant torus.
%For large helicity it might require a sufficiently high energy of the corresponding velocity field related to it by the ${\rm curl}^{-1}$  operator. 
Given two vorticity fields that differ by the number of hyperbolic zeros, one can arrange their arbitrary equal helicities, as well as equal large energies of the corresponding velocities, by using a pair of vortex plugs of opposite signs in each of the fields. 
Namely, for each of the vorticity fields in a small neighborhood in $M$ we insert two vortex plugs one after another and ``rotating" in the opposite directions. The first plug would  approximately ``add" vorticity helicity $W$ and  velocity energy $E$ (given by the bounded ${\rm curl}^{-1}$ operator), while the next one would add, respectively
$-W$ and $E$. Thus by considering linear combinations of those plugs with large positive coefficients 
one has the required two-parameter control of the helicity-energy integrals.
(For alternative arguments, also based on a combination of geometry and properties of  ${\rm curl}^{-1}$, cf. \cite{CT, KKP1}.) 

\medskip

$(ii)$ The easiest way to observe non-transitivity is to recall that the Euler equation in any dimension has such first integrals as generalized 
enstrophies  and helicities. They are defined with the help of 
the 1-form $u:=v^\flat$ corresponding to a velocity field $v$ by using the Riemannian metric on $M$. Namely, for an even-dimensional manifold $M^{2k}$ 
all moments of the $2k$-form $(du)^k$, i.e. all generalized 
enstrophies $I_m(v)=\int_M ((du)^k/\mu)^m\mu$, are first integrals of the Euler equation.
For an odd-dimensional $M^{2k+1}$ the generalized helicity $I(v)=\int_M u\wedge (du)^k$ is a first integral. 

Furthermore, a $C^1$ neighborhood of the velocity $v$ corresponds to a $C^0$ neighborhood of the vorticity 2-form $du$.
Therefore two sufficiently small   $C^1$ neighborhoods of the velocities $v_0$ and $v_1$ with different values of  generalized helicities (for odd dimension)
or generalized enstrophies (for even dimension) will never overlap during the Euler evolution. This implies non-transitivity for the problem with no extra constraints.
(Alternatively, one can use a local construction of invariants of solutions near singular points of their vorticity form described in Example 4.
Then initially non-intersecting $C^1$ neighborhoods of the vorticity with such singularities will never overlap.)
%This would allow one to find non-overlapping neighborhoods satisfying certain global constraints, such as equality of generalized helicity.}
%$\diamondsuit$
\proofend
\medskip

{\bf Example 2.} Here we describe a construction of a ``vortex plug", 
which is a local deformation (or surgery) of the vorticity field in 3D allowing one to change  its helicity by an arbitrary amount. 
Rectify the vorticity field in a neighborhood of a nonsingular point and 
consider a short invariant cylinder inside that neighborhood with cylindrical coordinates $(r, \theta, z)$ with the volume form $\mu:=r\,dr\,d\theta\,dz$, where $r\ge 0$.
We deform the vorticity field inside this invariant topological cylinder  in the following way to introduce two non-degenerate zeros (and a nondegenerate periodic orbit) while keeping the field divergence-free and without changing it on the boundary.  

First consider the 2D setting and a family of Hamiltonian fields with Hamiltonian functions $H(r, z)=r^2(r^2+z^2+a), \, a\in \mathbb R$ and non-standard symplectic structure $r\,dr\wedge dz$ degenerate along the line $r=0$. The corresponding Hamiltonian field has the form
$$
v_H:= -\frac 1r\frac{\partial H}{\partial z}\frac{\partial}{\partial r}+\frac 1r\frac{\partial H}{\partial r}\frac{\partial}{\partial z}
=-2rz\frac{\partial}{\partial r}+2(2r^2+z^2+a)\frac{\partial}{\partial z}\,.
$$
One can see that for a positive value of $a$  this Hamiltonian field 
is topologically equivalent to that for  the Hamiltonian $H=r^2$ near the origin. 
When the parameter $a$ changes from $0^+$ to $0^-$, two saddles and two centers are born for the corresponding Hamiltonian field $v_H$. 

Now we consider an axisymmetric  3D  analog of that Hamiltonian field by adding a special rotation about the $z$-axis, the field 
$v_f:=f(H){\partial}/{\partial \theta}$ which rotates each level of $H$ with its own speed. One can  see that the field $v_{H,f}=v_H+v_f$ 
is divergence-free, ${\rm div}_\mu v_{H,f}=0$ for any function $f$. Now, by setting, e.g. $a=-1$ and choosing an appropriate $f$ one obtains a field $v_{H,f}$ which has a unit sphere $\{r^2+z^2=1\}\subset \mathbb R^3$ as an invariant surface, two nondegenerate focus points  
$(r=0, z=\pm 1)$ on it, and a family of nested invariant tori inside, whose core is a nondegenerate periodic orbit born from the two centers.
By choosing an arbitrary speed of rotation in the $\theta$-direction by means of the function $f$, one can achieve an arbitrary helicity in that region, and it does not add to the helicity outside, since this invariant sphere is not linked with anything outside of it in $M$. We call it a ``vortex plug". 
By inserting such  plugs to the original vorticity field one can attain any prescribed helicity. One can also see that by changing the parameter $a$ one can include the plug construction into a deformation  from a field without zeros, rather than making it as a surgery.
 \medskip

{\bf Remark 3.}
An advantage of  the $C^1$-setting for vorticity is that  one has plenty of locally defined continuous Casimirs --
namely, multiplicators (i.e. eigenvalues of the linearization) of non-degenerate zeros of vorticity. 
Each nondegenerate zero of vorticity in 3D  with simple eigenvalues gives $2$ locally defined Casimirs (the eigenvalues at each singularity sum to zero because of the divergence-free condition). These Casimirs may replace the first integrals  of \cite{KKP1}
measuring the volume of invariant tori in a given isotopy class (or maybe other more subtle invariants for nonvanishing fields). In the $C^0 $-case one has only semicontinuous integer-valued Casimirs measuring the number of zeros 
or the like.  For instance, in \cite{CT} one studies vorticities of contact type, which must have no zeros, so the index argument provides a weaker requirement for a neighborhood to stay away from  fields with zeros, which are certain not to be  of contact type. In order to prove  c) one can consider fields $w_k$ that have at least $2k$ non-degenerate zeros, thus providing a countable number of neighbourhoods.

 Finally, note that for d) and e) one needs  to use more subtle arguments. For instance, 
 whenever one imposes additional constraints for a nonvanishing field, e.g. to stay in the 
same homotopy class, subtle invariants of contact homology are employed in \cite{CT}. The local non-mixing discussed
in \cite{KKP2} was based on a specific fibrated structure of 3D steady solutions. 
 
\medskip

{\bf Example 4.} Here we present a construction of $C^2$ local invariants for velocity in any dimension. 
For an odd-dimensional manifold $M^{2k+1}$ the Darboux theorem says that a generic (i.e. maximally non-integrable) 1-form $u$ 
locally has the form $u=dz+\sum x_idy_i$. Moreover, one can achieve the same form via volume-preserving transformations (see e.g. [1])
and thus this form has no local invariants. Consider now a special form $u=d(z^2)+\sum^k_{i=1} x^2_id(y^2_i)$.  We compute that 
$du= \sum^k_{i=1} d(x^2_i)\wedge d(y^2_i)$ and 
$$
u\wedge (du)^k= d(z^2) \bigwedge^k_{i=1} d(x^2_j)\wedge d(y_j^2) =2^{2k+1}x_1\dots y_kz\,  dx_1\wedge \dots \wedge dy_k  \wedge dz \,.
$$
In the latter expression the origin $x_1=...=y_k=z=0$ is a special point where the form vanishes. Moreover, it is stable: any $C^2$ perturbation of the original form $u$ will have the same type of degeneration for $u\wedge (du)^k$. And finally, any volume-preserving diffeomorphism of a neighborhood of that point cannot change the constant coefficient. In particular, e.g. the form $au$ is not volume-preserving diffeomorphic to $u$ in the vicinity of the origin for $a\not= 1$. On the other hand, adjusting $u$ outside of the origin, one can achieve any value of the total integral for the form  $u\wedge (du)^k$
i.e. of the total generalized helicity. 

Note also that in many cases one can extend this consideration to 
the $C^1$-setting for $u$. For instance, on a contact manifold $M^{2k+1}$ there is a 1-form $\tilde u$ such that $\tilde u\wedge (d\tilde u)^k$ has no zeros,
and hence its $C^1$ neighborhood can be separated from that of $u$ with zeros described above (this reminds the consideration  in \cite{CT} of 1-forms of contact type in 3D).

One can see that in 3D the degenerations of $du$ corresponds to zeros of the vorticity field $w={\rm curl}\,v$ defined by $i_w\mu =du$ for $u=v^\flat$. 
Then the local invariants of $u\wedge (du)^k$ at the origin 
are analogous to the multiplicators of the vorticity field $w$. 

\smallskip

Similarly, one constructs local invariants at degenerations of  $u=\sum^k_{i=1} x^2_id(y^2_i)$ on even-dimensional manifolds $M^{2k}$.
Note that in 2D they boil down to invariants of $du=4xy\,dx\wedge dy$, and hence to those of the Morse function $xy$, with respect to  transformations
 preserving  the area form $dx\wedge dy$. The front coefficient is an invariant delivered by Le Lemme de Morse Isochore \cite{Isochore}.

It is worth mentioning that without requirement of the volume preservation, local invariants of typical degenerations of closed 2-forms $du$  are subtle and  have been studied in \cite{Do}.
\medskip

{\bf Remark 5.}
To summarize, without imposing extra constraints, Problem \#31 from \cite{KMS} about non-mixing in low smoothness becomes rather straightforward, due to the existence of various Casimirs separating coadjoint orbits. Actually, this is a property of all Euler-Arnold equations for any Lie group: the existence of a Casimir which separates neighborhoods of coadjoint orbits implies  non-mixing of the corresponding equation in the dual of  its Lie algebra. This  property is based solely on the ``kinematics"  of the equation. So from this point of view any Euler-Arnold equation admitting Casimirs
(and in particular, the Euler fluid dynamics in any dimension) is non-mixing. 
The above  leads to the following natural rectification of the problem from \cite{KMS}: 
\smallskip

{\bf Problem 6.} Is the Euler equation non-mixing  within  the coadjoint orbits of  initial vorticity functions in 3D?
\smallskip

An answer to this question should involve the actual ``dynamics" of the Euler equation, similarly to the study of wandering orbits in 2D fluids in \cite{DEJ, Nad}, which is counterposed to finite dimensions with Poincar\'e's recurrence on all compact coadjoint orbits.  
 One may hope that  such tools as contact-type forms and KAM,  proposed in  \cite{CT, KKP1,KKP2},  might be useful for this updated problem or could  find other applications in hydrodynamics. 
\medskip

{\it Acknowledgements.} BK thanks Daniel Peralta-Salas and Theodore Drivas  for fruitful discussions. The research was supported by an NSERC Discovery Grant.

\end{document}